\documentclass[12pt,twoside]{amsart}
\usepackage{amssymb}
\usepackage{amscd}

\title{The indices of log canonical singularities}  
\author[Osamu Fujino]{Osamu Fujino*}
\thanks{*Research Fellow of the Japan Society for the Promotion of Science}
\address{Research Institute for Mathematical Sciences\\ 
Kyoto University, Kyoto 606-8502 Japan}
\email{fujino@kurims.kyoto-u.ac.jp}
\newcommand{\bQ}[0]{{\mathbb Q}}

\newcommand{\xDiff}[0]{{\operatorname{Diff}}}
\newcommand{\Supp}[0]{{\operatorname{Supp}}}

\newcommand{\xCLC}[0]{{\operatorname{CLC}}}
\newcommand{\xLLC}[0]{{\operatorname{LLC}}}
\newcommand{\xCLCC}[0]{{\operatorname{CLC^{c}}}}
\newcommand{\xLLCC}[0]{{\operatorname{LLC^{c}}}}

\newcommand{\xdiscrep}[0]{{\operatorname{discrep}}}
\newcommand{\xInd}[0]{{\operatorname{I\,}}}
\newcommand{\xAut}[0]{{\operatorname{Aut}}}
\newcommand{\xGL}[0]{{\operatorname{GL}}}
\newcommand{\xBim}[0]{{\operatorname{Bim}}}
\newcommand{\xBir}[0]{{\operatorname{Bir}}}

\newcommand{\ad}[0]{{admissible}}
\newcommand{\pad}[0]{{preadmissible}}
\newcommand{\bb}[0]{{B-birational}}
\newcommand{\bm}[0]{{B-bimeromorphic}}

\newcommand{\bp}[0]{{B-part}}
\newcommand{\bs}[0]{{B-smooth}}
\newcommand{\Spec}[0]{{\operatorname{Spec}}}
\newcommand{\xRes}[0]{{\operatorname{Res}}}
\newcommand{\xcenter}[0]{{\operatorname{Center}}}

\newcommand{\xlcm}[0]{{\operatorname{l.c.m.}}}
\newcommand{\wt}[0]{{\operatorname{wt}}}

\newtheorem{thm}{Theorem}[section]
\newtheorem{lem}[thm]{Lemma}
\newtheorem{cor}[thm]{Corollary}
\newtheorem{prop}[thm]{Proposition}

\newtheorem{conj}[thm]{Conjecture}
\newtheorem{ld}[thm]{Lemma-Definition}

\theoremstyle{definition}
\newtheorem{defn}[thm]{Definition}
\newtheorem{rem}[thm]{Remark}
\newtheorem*{ack}{Acknowledgements}       

\newtheorem{notation}[thm]{Notation}

\newtheorem{exmp}[thm]{Example}

\newtheorem*{notation2}{Notation}         

\theoremstyle{remark}

\begin{document}
\bibliographystyle{amsalpha+}
\maketitle

\abstract 
Let $(P\in X,\Delta)$ be a three dimensional log canonical pair 
such that $\Delta$ has only standard coefficients and 
$P$ is a center of log canonical singularities for $(X,\Delta)$. 
Then we get an effective bound of the indices of these pairs 
and actually determine all the possible indices. 
Furthermore, under certain assumptions including the log Minimal 
Model Program, an effective bound is also obtained in dimension $n\geq 4$. 
\endabstract

\setcounter{section}{-1}
\section{Introduction}\label{se0}
The main purpose of this paper is 
to investigate the indices of log canonical 
pairs. 
Let $(P\in X)$ be a log canonical singularity 
which is not log terminal. 
If $\dim X=2$, then the index is $1,2,3,4$, or $6$. 
This fact is well-known to specialists. 
Shihoko Ishii generalized this result to 
three dimensional isolated log canonical singularities which 
are not log terminal. 
More precisely, she proved that a positive integer $r$ is the index of 
such a singularity if and only if $\varphi (r)\leq 20$ and 
$r\ne 60$, where $\varphi$ is the Euler function (for related topics, 
see \cite{Sh2}). 
In this paper, we generalize it to higher dimensional 
(not necessarily isolated) log canonical singularities 
which are not log terminal. 
We note that if $(P\in X)$ is a log canonical singularity such that 
$P$ is not a center of log canonical singularities, 
then the index is not bounded (see Example (\ref{ex7})). 
So, we shall prove the following 
(for the precise statement, see Corollary (\ref{maincor}) and 
Remark (\ref{rem22})). 

\begin{thm}\label{thmi1}  
Let $(P\in X,\Delta)$ be a three dimensional log canonical pair such 
that 
$\Delta$ has only standard coefficients and 
$P$ is a center of log canonical singularities for the pair $(X,\Delta)$. 
Then the index of $(X,\Delta)$ at $P$ is bounded. 
More precisely, the positive integer $r$ is the index of such a pair if 
and only if $\varphi (r)\leq 20$ and 
$r\ne 60$. 
In particular, if there exists another center of log canonical singularities 
$W$ such that 
$P\subsetneq W$, then the index is $1,2,3,4$, or $6$. 
\end{thm}

This is related to (birational) automorphisms on 
$K3$ surfaces, Abelian surfaces and elliptic curves. 
Unfortunately, (birational) automorphisms on higher dimensional ``Calabi-Yau'' 
varieties are not well understood. 
If we can prove the conjectures about such automorphism groups 
(see Conjecture (\ref{conj})), then Theorem (\ref{thmi1}) is generalized 
to the following (for the precise and effective 
statement, see Theorem (\ref{mainthm})). 
Precisely, we prove Theorem (\ref {thmi2}) and get Theorem (\ref{thmi1}) 
as a corollary. 

\begin{thm}\label{thmi2}
Assume the log Minimal Model Program for dimension $\leq n$. 
Let $(P\in X,\Delta)$ be 
an $n$-dimensional log canonical pair such that 
$\Delta$ has only standard coefficients and 
$P$ is a center of log canonical 
singularities for the pair $(X,\Delta)$. 
If the conjectures $(F'_{n-1})$ and $(F_{l})$ hold 
true for $l\leq n-2$ {\em{(}}see Conjecture (\ref{conj}){\em{)}}, 
then the index of $(X,\Delta)$ at $P$ is bounded. 
\end{thm}

This theorem is an answer to \cite[4.16]{I}. 
We should mention 
that the idea of this paper is due to \cite{I2} and \cite{I}, 
and the proof relies on \cite{Fj} (see also \cite{Sh2}). 

We explain the contents of this paper. 
In Section \ref{se1}, we fix our notation and recall some definitions 
used in this paper, some of which were introduced in \cite{Fj}. 
In Section \ref{se2}, we generalize Shokurov's connectedness lemma. 
This section is a continuation of \cite [Section 2]{Fj}. 
Section \ref{se3} deals with birational automorphism groups and 
we collect some known results for low dimensional 
varieties. 
Section \ref{se4} 
is devoted to the proof of the main result, Theorem (\ref{mainthm}). 
In Section \ref{se5}, we collect some examples of log canonical singularities. 
Finally, in Section \ref{se6}, we explain how to translate statements 
on algebraic varieties into those on analytic spaces. 

\begin{notation2}
(1) We will make use of the standard notation and definitions as in 
\cite{KM}. 

(2) The log Minimal Model Program (log MMP, for short) means the log MMP 
for $\bQ$-factorial dlt pairs. 

(3) A {\em{variety}} means an algebraic variety over $\mathbb C$ and an 
{\em{analytic space}} a reduced complex analytic space. 
\end{notation2}

\begin{ack}
I would like to express my gratitude to 
Professor Shihoko Ishii, whose seminar talk at 
RIMS was the starting point of this paper. 
I am grateful to Professors 
Yoichi Miyaoka, Shigefumi Mori, and Noboru Nakayama 
for giving me many useful comments. 
I am also grateful to Mr.~Hiraku Kawanoue, 
especially, his advice about the desingularization theorem.  
\end{ack}

\section{Preliminaries}\label{se1}

In this section, we fix our notation and recall some definitions. 
For analytic spaces, we have to modify Definitions (\ref{def6}), 
(\ref{def3.1}) 
and Lemma (\ref{99}) (see Section {\ref{se6}). 

\begin{notation}\label{notation1}
Let $X$ be a normal variety over $\mathbb C$. 
The {\em{canonical divisor}} $K_X$ is defined so that 
its restriction to the regular part of $X$ is a divisor of a regular 
$n$-form. 
The reflexive sheaf of rank one $\omega _X:={\mathcal O}(K_X)$ 
corresponding to $K_X$ is called the {\em{canonical sheaf}}. 
Let $(P\in X,\Delta)$ be a germ of a normal variety with $\bQ$-divisor 
such that $K_X+\Delta$ is $\bQ$-Cartier. 
The {\em{index}} of $(X,\Delta)$ at $P$, denoted by $\xInd (P\in X,\Delta)$, 
is the smallest positive integer $r$ such that $r(K_X+\Delta)$ 
is Cartier at $P$. 
\end{notation}

The following is the definition of singularities of pairs. 
Note that the definitions in \cite {KMM} or \cite{KM} are 
slightly different from ours. 

\begin{defn}\label{def1}
Let $X$ be a normal variety and $D=\sum d_i D_i$ an effective 
$\bQ$-divisor such that $K_X+D$ is $\bQ$-Cartier. 
Let $f:Y\to X$ be a proper birational morphism. 
Then we can write 
$$
K_Y=f^{*}(K_X+D)+\sum a(E,X,D)E, 
$$ 
where the sum runs over all the distinct prime divisors $E\subset Y$, 
and $a(E,X,D)\in \bQ$. This $a(E,X,D)$ is called the 
{\em discrepancy} of $E$ with respect to $(X,D)$. 
We define 
$$
\xdiscrep (X,D):=\inf _{E}\{a(E,X,D)\ 
|\  E \text{ is exceptional over}\  X  \}.
$$
On the assumption that $0\leq d_i \leq 1$ 
for every $i$, we say that $(X,D)$ is 
$$
\begin{cases}
\text{terminal}\\
\text{canonical}\\
\text{klt}\\
\text{plt}\\
\text{lc}\\
\end{cases}
\quad {\text{if}} \quad \xdiscrep (X,D) 
 \quad
\begin{cases}
>0,\\
\geq 0,\\
>-1\quad {\text {and \quad $\llcorner D\lrcorner =0$,}}\\
>-1,\\
\geq -1.\\
\end{cases}
$$
Moreover, $(X,D)$ is {\em{divisorial log terminal}} (dlt, for short) if 
there exists a log resolution (see \cite [Notation 0.4 (10)]{KM}) with 
$a(E,X,D)>-1$ for every exceptional divisor $E$. 
Here klt (resp.~plt, lc) is short for {\em {Kawamata log terminal}} (resp.~
{\em{purely log terminal}}, {\em{log canonical}}). 
If $D=0$, then the notions klt, plt and dlt coincide and in this case we say 
that $X$ has {\em{log terminal}} (lt, for short) singularities. 

Let $S:=\{1-1/m\,|\,m\in {\mathbb{N}}\cup \{\infty\}\}$. 
We say that the divisor $D=\sum d_i D_i$ has {\em{only standard coefficients}} 
if $d_i\in S$ for every $i$ (cf.~\cite [1.3]{Sh2}). 
\end{defn}

In the following definition, we define the {\em 
{compact center of log canonical 
singularities}}. 

\begin{defn}[{cf.~\cite [Definition 1.3]{Ka}}]\label{def2}
A subvariety $W$ of $X$ is said to be a 
{\em center of log canonical singularities} 
for the pair $(X,D)$, if there exists a proper birational morphism 
from a normal variety $\mu:Y\to X$ and a prime divisor $E$ on 
$Y$ with the discrepancy coefficient $a(E,X,D)\leq -1$ such 
that $\mu(E)=W$. 

The set of all centers of log canonical singularities is 
denoted by $\xCLC (X,D)$. 
The union of all the subvarieties in $\xCLC (X,D)$ is denoted by 
$\xLLC (X,D)$ and called the 
{\em locus of log canonical singularities} for 
$(X,D)$. 
$\xLLC (X,D)$ is a closed subset of $X$. 

We denote the set of {\em {compact}} (with respect to 
the classical topology) elements in $\xCLC (X,D)$ by 
$\xCLCC (X,D)$. 
If $W\in \xCLCC(X,D)$, then $W$ is said to be a 
{\em {compact center of log canonical singularities}} 
for the pair $(X,D)$. 
It can be checked easily that $W\in \xCLCC(X,D)$ if and only if 
there exists a proper birational morphism 
from a normal variety $\mu:Y\to X$ and a compact prime divisor $E$ on 
$Y$ with the discrepancy coefficient $a(E,X,D)\leq -1$ such 
that $\mu(E)=W$. 
We define $\xLLCC (X,D)$ as the union of subvarieties in 
$\xCLCC (X,D)$. 
\end{defn}

\begin{rem}\label{dlttuika}
Let $(X,D)$ be a dlt pair. 
Then there exists a log resolution $f:Y\to X$ such that $f$ induces an 
isomorphism over every generic point of center of log canonical 
singularities for the pair $(X,D)$ and $a(E,X,D)>-1$ for every 
$f$-exceptional divisor $E$. 
This is obvious by the original definition of dlt (see \cite[1.1]{Sh}). 
See also \cite[Divisorial Log Terminal Theorem]{Sz}. 
\end{rem}

Definitions (\ref {def5}), (\ref{def6}), (\ref{def3.1}), (\ref{sdlt}), 
and (\ref{admissible}) are reformulations of the definitions of 
\cite{Fj}. 

\begin{defn}[{cf.~\cite[Definition 4.6]{Fj}}]\label{def5} 
Assume that $X$ is nonsingular and 
$\Supp D$ is a simple normal crossing divisor and 
$D=\sum_id_i D_i$ is a ${\bQ}$-divisor 
such that $d_i\leq 1$ ($d_i$ may be negative) for every $i$. 
In this case we say that $(X,D)$ is {\it {\bs}}. 

Let $(X,D)$ be dlt or {\bs}. 
We write $D=\sum_i d_i D_i$ 
such that $D_i$'s are distinct prime divisors. 
Then the {\it{\bp}} of $D$ 
is defined by $D^{B}:=\sum_{d_i=1}D_i$. 
We define the compact and non-compact {\bp} of $D$ as follows: 
$$ 
D^{c}:= \sum_{\scriptsize{\begin{array}{c} d_i=1\\  
D_i \text { is compact.}
\end{array}} }D_i\ ,\ \   
D^{nc}:=\sum_{\scriptsize{\begin{array}{c} d_i=1\\  
D_i \text { is non-compact.}
\end{array}} }D_i\, .  
$$

If $(X,D)$ is dlt or {\bs}, then a center of log canonical singularities 
is an 
irreducible component of an intersection of some {\bp} divisors. 
(See the Divisorial Log Terminal Theorem of \cite {Sz} and 
\cite [Section 2.3]{KM}.) 
When we consider a center of log canonical singularities 
$W$, we always consider the pair 
$(W,\Xi)$ such that $K_W+\Xi=(K_X+D)|_{W}$, 
where $\Xi$ is defined by repeatedly using the adjunction. 
Note that if $(X,D)$ is dlt (resp. ~{\bs}), 
then $(W,\Xi)$ is dlt (resp. ~{\bs}) by the adjunction. 

If $(X,D)$ is dlt or {\bs} and $W$ is a 
center of log canonical singularities for 
the pair $(X,D)$, then we write $(W,\Xi)\Subset (X,D)$. 
If there is no confusion, we write $W\Subset X$.  
\end{defn}

\begin{defn}[{cf.~\cite[Definition 1.5]{Fj}}]\label{def6}
Let $(X,D)$ and $(X',D')$ be normal varieties 
with $\bQ$-divisors such that $K_{X}+D$ and $K_{X'}+D'$
are $\bQ$-Cartier $\bQ$-divisors.

We say that 
$f: (X,D)\dashrightarrow (X',D')$ is a 
{\bb} map (resp.~morphism) if $f:X\dashrightarrow X'$ 
is a proper birational map
(resp.~morphism) and there exists a common resolution 
$\alpha:Y\to X$, $\beta:Y \to X'$ of 
$f:X\dashrightarrow X'$ such  that 
$\alpha ^{*}(K_{X}+D)=\beta ^{*}(K_{X'}+D')$.  
\end{defn}

The following lemma, which is a corollary of \cite [Claims ($A_n$), ($B_n$)]
{Fj}, is very useful.  

\begin{lem}\label{99}
Let $(X_i,D_i)$ be dlt or {\bs} for $i=1,2$. 
Let $f:(X_1,D_1)\dashrightarrow (X_2,D_2)$ be a {\bb} map and 
$(W_1,\Xi_1)\Subset (X_1,D_1)$ a minimal {\em{(}}with 
respect to $\Subset${\em{)}} center of log canonical singularities. 
Then there exists a minimal center of log canonical singularities 
$(W_2,\Xi_2)\Subset (X_2,D_2)$ such that 
$(W_1,\Xi_1)$ and $(W_2,\Xi_2)$ are $B$-birationally equivalent to 
each other. 
\end{lem}
\proof 
By Remark (\ref{dlttuika}), we may assume that $(X_i, D_i)$ is {\bs} for 
$i=1,2$. 
Let $g_i:(Y,E)\to (X_i,D_1)$ be a common resolution of $f$. 
We note that 
$g_1^{*}(K_{X_1}+D_1)=K_Y+E=g_2^{*}(K_{X_2}+D_2)$. 
By \cite [Claim ($A_n$)]{Fj}, we can take $(W,\Xi)\Subset (Y,E)$ such 
that 
$g_1|_{W}:(W,\Xi)\to (W_1,\Xi_1)$ is {\bb}. 
It is obvious that $(W,\Xi)\Subset (Y,E)$ is a minimal center of log 
canonical singularities for the pair $(Y,E)$ since 
$(W,\Xi)$ and $(W_1,\Xi_1)$ have the same discrepancies. 
By applying \cite[Claim ($B_n$)]{Fj} (see also Section \ref{se6}) 
to $g_2:(Y,E)\to (X_2,D_2)$, $(W,\Xi)\to (W_2,\Xi_2)$ is {\bb}, 
where $W_2=g_2(W)$ and $(W_2,\Xi_2)\Subset (X_2,D_2)$. 
It is obvious that $(W_2,\Xi_2)$ is a minimal center of log canonical 
singularities for the pair $(X_2.D_2)$ by \cite [Claim ($A_n$)]{Fj}.    
\endproof

\begin{defn}[{cf.~\cite[Definition 3.1]{Fj}}]\label{def3.1}
Let $(X,D)$ be a pair of a normal variety and a ${\bQ}$-divisor 
such that $K_X+D$ is ${\bQ}$-Cartier. We define 
\begin{eqnarray*}
\xBir(X,D)&:=& \{\sigma :(X,D)\dashrightarrow (X,D)\ |\ \sigma 
\text { is a {\bb} map}\},\\
\xAut(X,D)&:=& \{\sigma :X\to X\ |\ 
\sigma \text { is an automorphism and } 
\sigma^{*}D=D\}.
\end{eqnarray*}
Since $\xBir(X,D)$ acts on 
$H^{0}(X,{\mathcal O}_{X}(m(K_X+D)))$ for every integer $m$ 
such that $m(K_X+D)$ is a Cartier divisor,
we can define B-pluricanonical representations 
$\rho_{m}:\xBir(X,D)\to \xGL(H^{0}(X,m(K_X+D)))$.
\end{defn}

The following is the definition of semi divisorial log terminal 
pairs. 
The notion of semi divisorial log terminal is much better than 
that of semi log canonical for the inductive treatment (see \cite{Fj}). 

\begin{defn}[{cf.~\cite[Definition 1.1]{Fj}}]\label{sdlt}
Let $X$ be a reduced algebraic scheme, which satisfies $S_2$ condition. 
We assume that it is
pure $n$-dimensional and 
normal crossing in codimension 1.
Let $\Delta$ be an effective $\bQ$-Weil divisor on $X$ 
(cf.~\cite[16.2 Definition]{FA}) such that 
$K_X+\Delta$ is $\bQ$-Cartier.

Let $X=\cup X_i$ be a decomposition into irreducible components,
and $\mu : X':= \amalg X'_i \to X=\cup X_i$ the normalization.
A $\bQ$-divisor $\Theta$ on $X'$ is defined by
$K_{X'} +\Theta:= \mu^*(K_X+\Delta)$ and a $\bQ$-divisor
$\Theta_i$ on $X'_i$ by $\Theta_i:=\Theta|_{X'_i}$.

We say that $(X,\Delta)$ is a {\it semi divisorial log terminal $n$-fold}
(an sdlt $n$-fold, for short)
if $X_i$ is normal, that is, $X'_i$ is isomorphic to $X_i$, 
and $(X',\Theta)$ is dlt.
\end{defn}

\begin{defn}[{cf.~\cite[Definition 4.1]{Fj}}]\label{admissible}
Let $(X,\Delta)$ be a proper sdlt $n$-fold
and $m$ a divisible integer.
We define {\ad} and {\pad} sections inductively on dimension.
\begin{itemize}
\item $s\in H^{0}(X,m(K_X+\Delta))$ is {\it{\pad}} 
if the restriction $s|_{(\amalg _{i}{\llcorner\Theta_{i}\lrcorner})}\in 
H^{0}(\amalg _{i} {\llcorner\Theta_{i}\lrcorner},m(K_{X'}+\Theta)|_
{(\amalg _{i}{\llcorner\Theta_{i}\lrcorner})})$
is {\ad}.
\item $s\in H^{0}(X,m(K_X+\Delta))$ is {\it{\ad}}
if $s$ is {\pad} and $g^{*}(s|_{X_j})=s|_{X_i}$ for 
every {\bb} map $g:(X_i,\Theta_i)\dashrightarrow (X_j,\Theta_j)$ for 
every $i,j$.
\end{itemize}
Note that if $s\in H^{0}(X,m(K_X+\Delta))$ 
is {\ad}, then the restriction $s|_{X_i}$ is 
$\xBir(X_i,\Theta_i)$-invariant for every $i$.
\end{defn}

The next lemma-definition is frequently used in Section \ref{se4}. 

\begin{ld}\label{ld8}
Let $(o\in X,\Theta)$ be a pointed variety with $\bQ$-divisor $\Theta$ 
such that $K_X+\Theta$ is $\bQ$-Cartier. 
Then there exists a resolution $f:(Y,\Xi)\to (o\in X,\Theta)$ such that 
\begin{enumerate}
\item [(1)] $f$ is projective, 
\item [(2)] $K_Y+\Xi=f^{*}(K_X+\Theta)$,
\item [(3)] $f^{-1}(o)$ is a simple normal crossing divisor in $Y$, 
\item [(4)] $f^{-1}(o)\cup \Xi$ is also a simple normal crossing 
divisor in $Y$.  
\end{enumerate}
We say that $f:(Y,\Xi)\to (o\in X,\Theta)$ 
is a very good resolution of $(o\in X,\Theta)$. 
\end{ld}

\proof
Let $f_1:Y_1\to X$ be any log resolution such that $f_1$ is projective. 
Apply the embedded resolution to $f_{1}^{-1}(o)_{\text {red}} \subset Y_1$. 
Then we get a sequence of blowing-ups $g:Y_2\to Y_1$. 
By Hironaka's theorem (see \cite {Hi} or \cite{BM}), the proper transform of 
$f_{1}^{-1}(o)_{\text {red}}$ by $g$ in $Y_2$ is smooth and the exceptional 
locus of $g$ is a simple normal crossing divisor, which intersects 
the proper transform of $f_{1}^{-1}(o)_{\text {red}}$ transversally. 
We note that $g$ is an isomorphism over $Y_1\setminus 
f_{1}^{-1}(o)_{\text {red}}$. 
Therefore, all the components of $(g\circ f_1)^{-1} (o)_{\text {red}}$ of 
codimension $\geq 2$ in $Y_2$ is smooth. 
By blowing up these components, we obtain  
$f_3:Y_3\to Y_2\to Y_1\to X$ and 
$\Xi_3$ such that $K_{Y_3}+{\Xi_3}={f_3}^{*}(K_X+\Theta)$ 
and $f_3^{-1}(o)$ is pure codimension $1$ in $Y_3$. 
We note that $Y_3$ is smooth. 
By applying the embedded resolution to 
$(f_{3}^{-1}(o)\cup \Xi_3)_{\text {red}}$, we get 
$f:Y\to Y_3\to Y_2\to Y_1\to X$ and $\Xi$ such that 
$K_Y+\Xi=f^{*}(K_X+\Theta)$. 
By the construction of $f$, 
$f$ is projective and $f:(Y,\Xi)\to (o\in X,\Delta)$ 
satisfies the conditions (3) and (4). 
\endproof 

The following lemma-definition follows from \cite [8.2.2 Lemma, 
17.10 Theorem]{FA}. 

\begin{ld}[$\bQ$-factorial dlt model]\label{dltmodel}
Assume that the log MMP holds in dimension $n$. 
Let $(X,\Theta)$ be an lc $n$-fold. 
Then there exists a projective birational morphism 
$f:(Y,\Xi)\to (X,\Theta)$ from a $\bQ$-factorial dlt pair $(Y,\Xi)$ 
such that $K_Y+\Xi=f^{*}(K_X+\Theta)$. 
Furthermore, if $(X,\Theta)$ is dlt, then we may take $f$ a small 
projective morphism. 
We say that $(Y,\Xi)$ is a $\bQ$-factorial dlt model of $(X,\Theta)$. 
\end{ld}

The following lemma is a special case of \cite[17.10 Theorem]{FA}. 
We use this in Lemma (\ref{taisetu}). 

\begin{lem}\label{taisetu2}
Assume the log MMP in dimension $n$. 
Let $(o\in Y,D)$ be a germ of a log canonical singularity such that 
$o\in \xCLC(Y,D)$. Without loss of generality, we may assume that 
$\xLLCC(Y,D)=o$. 
Let $h:(V,F)\to (o\in Y,D)$ be a very good resolution. 
Let 
$$
V = V^0 \stackrel{p_1}{\dashrightarrow} V^1
 \stackrel{p_2}{\dashrightarrow} \cdots
 \stackrel{p_i}{\dashrightarrow} V^i
 \stackrel{p_{i+1}}{\dashrightarrow} V^{i+1}
 \stackrel{p_{i+2}}{\dashrightarrow}
 \cdots \stackrel{p_{l-1}}{\dashrightarrow}
 V^{l-1} \stackrel{p_l}{\dashrightarrow} V^l = Z
$$
be the $(K_V+G)$-log MMP over $Y$, where $F:=G-H$ such that $G$ and $H$ 
are both effective $\bQ$-divisors without common irreducible components. 
We denote $F_{0}=F, G_{0}=G, H_{0}=H$, and 
$F_{i}=p_{i*}F_{i-1}, G_{i}=p_{i*}G_{i-1}, 
H_{i}=p_{i*}H_{i-1}$, for every $i$ and $F_{l}=E$. 
Then we obtain that 
$H_{l}=0$ and $f:(Z,E)\to (Y,D)$ is a $\bQ$-factorial dlt model, 
and $g:=p_{l}\circ p_{l-1}\circ\cdots\circ p_{1}$ induces an isomorphism 
at every generic point of center of log 
canonical singularities for the pair $(V,F)$. 
Furthermore, $\xLLCC(Z,E)=E^{c}$. 
\end{lem} 
\proof 
Since $h$ is a very good resolution, $\xLLCC(V,F)=F^{c}$ and 
$H$ contains no centers of log canonical singularities 
for the pair $(V,F)$. 
We note that 
$\xLLC(V,F)=\xLLC(V,G)$. 
By induction on $i$, 
we assume that $\xLLCC(V^{i},G_{i})=F_{i}^{c}=G_{i}^{c}$ and 
$g_i:= p_{i}\circ p_{i-1}\circ\cdots\circ p_{1}$ induces an 
isomorphism at every 
generic point of center of log 
canonical singularities for the pair $(V,F)$, and 
$H_{i}$ contains no centers of log canonical singularities 
for the pair $(V^{i},G_{i})$. 

If $p_{i+1}$ is a divisorial contraction, then $p_{i+1}$ contracts 
an irreducible component of $H_{i}$. 
Thus $\xLLCC(V^{i+1},G_{i+1})=F_{i+1}^{c}=G_{i+1}^{c}$. 
It is obvious that $H_{i+1}$ contains no 
centers of log canonical singularities 
for the pair $(V^{i+1},G_{i+1})$ and $g_{i+1}$ induces an 
isomorphism at every 
generic point of center of log 
canonical singularities for the pair $(V,F)$.  

If $p_{i+1}$ is a flip, then the flipping locus is included in 
$H_{i}$. In particular, every 
divisor whose center is in the flipping locus 
has discrepancy $>-1$ with respect to $K_{V^{i}}+G_{i}$. 
After the flip $p_{i+1}$, the discrepancies do not decrease. 
Therefore, $\xLLCC(V^{i+1},G_{i+1})=F_{i+1}^{c}=G_{i+1}^{c}$. 
Of course, $p_{i+1}$ is an isomorphism at every generic point of 
center of log canonical singularities for the pair $(V^{i},G_{i})$. 
In particular, $H_{i+1}$ contains no 
centers of log canonical singularities 
for the pair $(V^{i+1},G_{i+1})$  

Thus we get the result by the induction on $i$. 
We note that $H_{l}=0$ and $G_{l}=F_{l}=E$.  
\endproof

\begin{rem}\label{rr}
Let $h:(V,F)\to (o\in Y,\Delta)$ be a log resolution such that 
$K_{V}+F=h^{*}(K_{Y}+D)$. Assume that $h^{-1}(o)$ is not 
pure codimension $1$. 
Then $\xLLCC(V,F)$ does not necessarily coincide with $F^{c}$. 
We note that if there exist two irreducible components $F', F''$ 
of $F^{nc}$ such that $F'\cap F''\in \xCLCC(V,F)$, 
then $ F^{c}\subsetneq \xLLCC(V,F)$. 
This is why we need the notion of the very good resolution. 
\end{rem}

\section{Connectedness Lemmas}\label{se2} 

In this section, we treat connectedness lemmas. 
They play important roles in Section \ref{se4}. 
The results are stated for algebraic varieties. However, 
by the same argument, we can generalize them for analytic spaces 
which are projective over analytic germs (see Section {\ref{se6}). 
The following lemma is well-known (for the proof, see 
\cite [17.4 Theorem]{FA} or \cite[Theorem 1.4]{Ka}).   

\begin{lem}[{Connectedness Lemma, cf.~\cite [5.7]{Sh}, \cite[17.4]{FA}}] 
\label{lem3} 
Let $X$ and $Y$ be normal varieties and $f:X\to Y$ 
be a proper surjective morphism with connected fibers. 
Let $D$ be a ${\bQ}$-divisor on $X$ such that $K_X+D$ is ${\bQ}$-Cartier. 
Write $D=\sum d_i D_i$, where $D_i$ is an irreducible component of $D$ 
for every $i$. 
Assume the following conditions:  
\begin{enumerate}
\item [(1)] if $d_i<0$, then $f(D_i)$ has codimension at least two 
in $Y$, and 
\item [(2)] $-(K_X+\Delta)$ is $f$-nef and $f$-big.
\end{enumerate}
Then $\xLLC (X,D)\cap f^{-1}(y)$ is connected for every point $y\in Y$. 
\end{lem}

The next proposition is also well-known. 

\begin{prop}[{cf. \cite [6.9]{Sh} and \cite[12.3.2]{FA}}]\label{prop10}
Let $(S,\Theta)$ be a dlt surface. 
Let $f:S\to R$ be a proper surjective morphism onto a smooth curve $R$ 
with connected fibers. 
Assume that $K_S+\Theta$ is numerically $f$-trivial. 
Then $\llcorner \Theta\lrcorner \cap f^{-1}(r)$ has at most 
two connected components for every $r\in R$. 
Moreover, if $\llcorner \Theta\lrcorner \cap f^{-1}(o)$ has 
exactly two connected components for $o\in R$, 
then, in a neighborhood of $f^{-1}(o)$, $(S,\Theta)$ is 
plt and $\llcorner \Theta\lrcorner$ has no vertical component with respect to 
$f$. 
\end{prop}

By using Proposition (\ref {prop10}), we can prove the next corollary 
easily. 

\begin{cor}\label{prop11}
Let $(T,\Xi)$ be a {\bs} surface and $(S,\Theta)$ be a dlt surface. 
Let $g:(T,\Xi)\dashrightarrow (S,\Theta)$ be a {\bb} map 
and $h:T\to R$ and $f:S\to R$ be proper surjective morphisms onto a smooth 
curve $R$ with connected fibers such that $h=f\circ g$. 
Assume that there exists an irreducible component $C$ of 
$\Xi ^{B}$ such that 
$C\subset h^{-1}(r)$ for some $r\in R$. 
Then $\Xi ^{B}\cap h^{-1}(r)=\xLLC(T,\Xi)\cap h^{-1}(r)$ is connected. 
\end{cor}

\proof
By applying the elimination of indeterminacy, 
we may assume that $g$ is a {\bb} morphism. 
Apply Proposition (\ref{prop10}) to $f:S\to R$ and Lemma (\ref {lem3}) 
to $g:T\to S$. Thus we get $\Xi ^{B}\cap h^{-1}(r)$ is connected.  
\endproof

On the assumption that the log MMP holds in dimension $n$, 
we get higher dimensional generalizations of Proposition 
(\ref {prop10}) and Corollary (\ref{prop11}). 
Proposition (\ref{conn}) is a special case of \cite[Proposition 2.1 (0)]{Fj}. 

\begin{prop}\label{conn}
Assume that the log MMP holds in dimension $n$. 
Let $(X,\Theta)$ be a proper connected 
dlt $n$-fold and $K_X+\Theta \equiv 0$. 
Assume that $K_X+\Theta$ is Cartier. In particular, $\Theta$ is an 
integral divisor. 
Then  one of the following holds: 
\begin{enumerate}
\item[(1)] $\Theta$ is connected, 
\item[(2)] $\Theta=\Theta_1\amalg \Theta_2$, where 
$\Theta_i$ is connected and irreducible for $i=1,2$. 
In particular, $(X,\Theta)$ is plt. 
We note that $(X,\Theta)$ is canonical since $K_X+\Theta$ is Cartier and plt.  
Furthermore, there exists a {\bb} map $(\Theta_1,0)\dashrightarrow 
(\Theta_2,0)$. 
\end{enumerate}
\end{prop}

\proof
By replacing $(X,\Theta)$ with its $\bQ$-factorial 
dlt model, we may assume that $X$ is $\bQ$-factorial. 
By \cite [Proposition 2.1 (0) and Remark 2.2 (2)]{Fj}, 
$\llcorner\Theta\lrcorner=\Theta$ has at most 
two connected components. 
If $\Theta$ is connected, then (1) holds. 
So, we may assume that $\Theta$ has two connected components. 
By the proof of \cite[Proposition 2.1]{Fj}, 
there exists a sequence of 
$(K+\Theta-\varepsilon \Theta)$-flips 
and divisorial contractions $p:X\dashrightarrow Z$ and 
$(K_Z+p_*\Theta-\varepsilon 
p_*\Theta)$-Fano contraction to $(n-1)$-dimensional 
lc pair $(V,P)$, denoted by $u:Z\to V$, where $P$ is the divisor 
such that $K_Z+p_*\Theta=u^*(K_V+P)$ (see \cite [Lemma 2.3]{Fj}). 
We denote $p_*\Theta=\Theta'_1\amalg \Theta'_2$. 
Then $(\Theta'_i, \xDiff (p_*\Theta -\Theta'_i))\simeq (V,P)$ 
for $i=1,2$. 
It is because $\Theta'_i\simeq V$ by Zariski's Main Theorem. 
Since $K_Z+p_*\Theta$ is Cartier and $p_*\Theta=\llcorner p_*\Theta\lrcorner$, 
and $(Z,p_*\Theta)$ is lc, 
$\xDiff (p_*\Theta-\Theta'_i)$ is an integral divisor. 
Therefore, the divisor $P$ is also integral. 
Since $p_*\Theta$ has no vertical component with respect to $u$, 
we have $P=0$. 
By \cite[Appendix]{N2} or \cite[Corollary 4.5]{Fj2}, $(V,0)$ is lt. 
Therefore, $(\Theta'_i, \xDiff (p_*\Theta -\Theta'_i))=(\Theta'_i,0)$ is 
lt for $i=1,2$. 
Then $(Z,p_*\Theta)$ is plt in a neighborhood of $p_*\Theta$. 
Thus, we get (2) (see \cite[Proposition 2.1, Lemma 2.3, and Lemma 2.4]{Fj}).  
\endproof

The following proposition is a higher dimensional analogue of 
Corollary (\ref{prop11}). 
The proof is similar to that of Proposition (\ref{conn}).  

\begin{prop}\label{conn2}
Assume that the log MMP holds in dimension $n$. 
Let $(X,\Theta)$ be a dlt $n$-fold and 
$f:X\to R$ be a proper surjective morphism onto a normal variety $R$ 
with connected fibers. 
Assume that $\dim R\geq 1$, $K_X+\Theta$ is Cartier, and 
$K_X+\Theta$ is numerically $f$-trivial. 
Let $o\in R$ be a closed point. 
Assume that there exists an irreducible component $\Theta_o$ of $\Theta$ 
such that $f(\Theta_o)=o$. 
Then $f^{-1}(o)\cap \Theta$ is connected. 
\end{prop}
\proof
As in the proof of Proposition (\ref {conn}), 
We may assume that $X$ is $\bQ$-factorial. 
By \cite [Proposition 2.1 (1)]{Fj}, 
it is enough to think about the case where there exists a sequence of 
$(K+\Theta-\varepsilon \Theta)$-flips 
and divisorial contractions $p:X\dashrightarrow Z$ 
over $R$ and 
$(K_Z+p_*\Theta-\varepsilon p_*\Theta)$-Fano 
contraction to $(n-1)$-dimensional 
lc pair $(V,P)$ over $R$, denoted by $u:Z\to V$. 
In this case, $\xcenter _Z\Theta _o$ is in $p_*\Theta \cap h^{-1}(o)$, 
where $h:Z\to V\to R$. 
If $p_*\Theta \cap h^{-1}(o)$ is connected, then we get the result 
(see \cite [Lemma 2.4]{Fj}). 
So, we may assume that $p_*\Theta \cap h^{-1}(o)$ is not connected. 
By shrinking $R$ to a small analytic neighborhood of $o\in R$, 
we may assume that $p_*\Theta=\Theta'_1\amalg \Theta'_2$ in a neighborhood of 
$h^{-1}(o)$. 
Without loss of generality, we may assume that  $\xcenter _Z\Theta _o\subset 
\Theta'_1$. 
In a neighborhood of $h^{-1}(o)$, $u|_{\Theta'_i}:\Theta'_i\to V$ 
is finite for $i=1,2$. 
It is because, if $u|_{\Theta'_i}:\Theta'_i\to V$ is not finite, 
then $\Theta'_1\cap\Theta'_2\ne\emptyset$ since 
$\Theta'_{3-i}$ is $u$-ample. 
By using Zariski's Main Theorem (see, for example, 
\cite[Theorem 1.11]{Ue}), we get $\Theta'_i\simeq V$ and 
$(\Theta'_i, \xDiff (p_*\Theta -\Theta'_i))\simeq (V,P)$ in a neighborhood of 
$h^{-1}(o)$ for $i=1,2$. 
Since $p_*\Theta$ has no $u$-vertical component in a neighborhood of 
$h^{-1}(o)$, we get $P=0$. 
Note that $P$ is integral (see the proof of Proposition (\ref{conn})). 
Then $(\Theta'_i, \xDiff (p_*\Theta -\Theta'_i))\simeq (V,0)$ in a 
neighborhood of $h^{-1}(o)$. 
Therefore, $(\Theta'_i,\xDiff (p_*\Theta -\Theta'_i))$ is lt 
in a neighborhood of $h^{-1}(o)$. 
So, $(Z,p_*\Theta)$ is plt in a neighborhood of $h^{-1}(o)$. 
This contradicts the assumption that $\xcenter _Z\Theta _o\subset 
\Theta'_1$. 
So, $p_*\Theta$ is connected in a neighborhood of $h^{-1}(o)$. 
Then we get the result.  
\endproof

The next proposition plays an essential role in the proof of 
the main theorem 
(see Proposition (\ref{prop4.2})). 

\begin{prop}\label{keyprop}
Assume that the log MMP holds in dimension $\leq n$. 
Let $(T,\Xi)$ be a $B$-smooth $n$-fold and 
$(S,\Theta)$ a dlt $n$-fold. 
Let $g:(T,\Xi)\dashrightarrow (S,\Theta)$ be a {\bb} map and 
$h: (T,\Xi)\to (o\in R)$ and $f:(S,\Theta) \to (o\in R)$ 
be proper surjective morphisms 
with connected fibers onto a 
germ $(o\in R)$ of a normal variety $R$ with $\dim R\geq 1$. 
Assume the following conditions: 
\begin{enumerate}
\item[(1)] $K_S+\Theta$ and $K_T+\Xi$ are Cartier divisors, 
\item[(2)] $K_S+\Theta\equiv _f 0$ and $K_T+\Xi\equiv _h 0$, 
\item[(3)] $h^{-1}(o)$ and  $h^{-1}(o)\cup \Xi$ are 
simple normal crossing divisors, 
\item[(4)] $g$ induces an isomorphism at every generic point of 
center of log canonical singularities for the pair $(T,\Xi)$, 
\item[(5)] there exists an irreducible component $\Xi_o$ of $\Xi^{B}$ 
such that $h(\Xi_o)=o\in R$.  
\end{enumerate}
Then $\xLLCC(T,\Xi)=\Xi^{c}$ is connected. 
Furthermore, if $\Xi'\subset \Xi^{B}$ is an irreducible component such that 
$h^{-1}(o)\cap \Xi'\ne \emptyset$, 
then $\Xi'\cap \Xi^{c}\ne\emptyset$. 
\end{prop}
\proof
First, if $\dim T=\dim S=2$, then this proposition is true by 
Lemma (\ref {lem3}) and Corollary (\ref {prop11}). 

Next, apply the elimination of indeterminacy. 
We may assume that $g$ is a morphism. 
By using 
Proposition (\ref{conn2}) 
and applying Lemma (\ref{lem3}) to the morphism $g$, 
we get $\xLLC(T,\Xi)\cap h^{-1}(o)=\Xi^{B}\cap h^{-1}(o)$ is connected. 

Finally, we go back to the original {\bb} map $g$. 
Let $D$ be the maximum connected component of $\Xi^{c}$ such that 
$\Xi_o\subset D$. 
For the proof of this proposition, 
it is enough to exclude the following situation; 
\begin{enumerate}
\item[$(\clubsuit)$] there exist irreducible divisors $\Xi_1$ and $\Xi_2$ 
which satisfy the following conditions:
\begin{enumerate}
\item[(i)] $\Xi_1$ and $\Xi_2$ are irreducible components of $\Xi^{B}$, 
\item[(ii)] there exists a connected component $C$ of 
$\Xi_1\cap h^{-1}(o)$ such that $C\cap D\ne \emptyset$ and 
$C\cap (\Xi_2\cap h^{-1}(o))\ne \emptyset$, 
\item[(iii)] $(\Xi_2\cap h^{-1}(o))\cap D =\emptyset$.
\end{enumerate}
\end{enumerate}
Note that $D\subset h^{-1}(o)$. 
We use the induction on $n$ to exclude $(\clubsuit)$. 
By the above argument, when $\dim T=\dim S=2$, this proposition is 
true. 
Assume that this proposition is true in dimension $n-1$. 
If $h(\Xi_1)=o$, then $\Xi_1\subset \Xi ^{c}$. 
This contradicts the definition of $D$. 
So we get $\Xi_1\subset \Xi ^{nc}$. 
Let $(S',\Theta')$ be a proper transform of 
$(\Xi_1,(\Xi-\Xi_1)|_{\Xi_1})$, which can be taken by the condition 
(4), and 
$h': (\Xi_1,(\Xi-\Xi_1)|_{\Xi_1})\to R'$ be the Stein factorization of 
$h:\Xi_1\to R$. 
Since $S'$ is normal, there exists $f':S'\to R'$. 
We define $o':=h'(C)$. 
Apply the hypothesis of the induction to 
$(\Xi_1,(\Xi-\Xi_1)|_{\Xi_1}),(S',\Theta')$, and $(o'\in R')$. 
We note that the conditions (1), (2), and (4) are satisfied by 
the adjunction. 
The condition (3) is true since $\Xi_1$ is an irreducible 
component of $\Xi$ and (5) is also true since $C\cap D\ne \emptyset$. 
Therefore, we obtain that, in the fiber $(h')^{-1}(o')$, $\Xi_1\cap D$ and 
$\Xi_1\cap \Xi_2$ are connected by $((\Xi-\Xi_1)|_{\Xi_1})^{c}$. 
Thus there exists $\Xi_3\subset \Xi^{c}$ such that $\Xi_3\cap D\ne 
\emptyset$, and $\Xi_3\not\subset D$. 
This contradicts the definition of $D$. 
We note the condition (3) and $\xLLCC(T,\Xi)=\Xi^{c}$. 
\endproof

\section{Finiteness and Boundedness}\label{se3}

In this section, we investigate the birational automorphism groups. 
First, we prove the following proposition, which is an easy consequence of 
\cite [Theorem 14.10]{Ue} and a special case of \cite[Conjecture 3.2]{Fj}.  

\begin{prop}\label{prop24}
Let $(S,0)$ be a normal $n$-fold with only canonical singularities 
such that $K_S\sim 0$. Then the image $\rho_m (\xBir (S,0))$ is finite, 
where $\rho_m: \xBir (S,0)\to \xGL (H^0(S,mK_S))$. 
\end{prop}

\proof
Let $f:T\to S$ be a resolution. 
Then $K_T=f^*K_S+E$, where $E$ is an effective Cartier divisor. 
It is obvious that $\sigma \in \xBir (S,0)$ induces $f^{-1}\circ 
\sigma \circ f\in \xBir (T,-E)$. 
We denote 
$$
\xBir (T)=\{ {\text{all birational maps of $Y$ onto itself}}\}.
$$ 
Then $\xBir (T,-E)\subset \xBir (T)$. 
The image of $\xBir (T)\to \xGL (H^{0}(T,mK_T))$ is 
finite 
by \cite[Theorem 14.10]{Ue}, 
Thus the image of $\xBir (T,-E)\to \xGL (H^{0}(T,m(K_T-E)))$ is also finite. 
Note that $H^{0}(T,m(K_T-E))=H^{0}(T,mK_T)$. 
Therefore $\rho_m (\xBir (S,0))$ is finite. 
\endproof

By Proposition (\ref{prop24}), we get the finiteness of $B$-pluricanonical 
representations. However, for the main theorem, we need the stronger results. 
So, we write down the required conjectures. 

\begin{conj}[\,Boundedness of $B$-canonical\, representations\,]
\label{conj}
The following conjectures $(F_l)$ and $(F'_l)$ are used in the main 
results.  
\begin{enumerate}
\item[($F_l$)] There exists a positive integer $B_l$ such that 
$|\rho_1(\xBir (S,0))|\leq B_l$ 
for every $l$-dimensional variety $S$ with only canonical singularities 
such that $K_S\sim 0$.  
\item[($F'_l$)] There exists a positive integer $B'_l$ such that 
the order of $\rho_1(g)$ in $\xGL (H^{0}(S,K_S))$ is bounded above by 
$B'_l$ 
for every $l$-dimensional variety $S$ with 
only canonical singularities such that 
$K_S\sim 0$ and for every 
$g\in \xAut (S,0)=\xAut (S)$ such that $g$ has a finite order.  
\end{enumerate}
\end{conj}

For low dimensional varieties, we know the details of 
canonical representations. 
We list the results needed in this paper for reader's convenience. 
This is \cite[Proposition 4.8, Proposition 4.9]{I}.   

\begin{prop}\label{prop20}
{\em{(1)}} For an arbitrary elliptic curve $C$, denote the order 
$|\rho_1(\xAut(C))|$ by $r$, where 
$\rho_1:\xAut (C)\to \xGL(H^{0}(C,K_{C}))$. 
Then $\varphi (r)\leq 2$, which means $r=1,2,3,4,6$ 
{\em{(}}see, for example, \cite[Chapter IV Corollary 4.7]{H}{\em{)}}. 
Here $\varphi$ is the Euler function. 

{\em{(2)}} 
{\em{(}}\cite[10.1.2]{Ni}{\em{)}} For an arbitrary $K3$ surface $X$, 
denote the order $|\rho_1(\xAut(X))|$ by $r$, 
where $\rho_1:\xAut (X)\to \xGL(H^{0}(X,K_{X}))$ is 
the induced representation. 
then $\varphi (r)\leq 20$, in particular $r\leq 66$. 
Machida and Oguiso proved that there are no $K3$ surfaces which satisfy 
$r=60$ in \cite{MO}. 

{\em{(3)}} 
{\em{(}}\cite [3.2]{Fk}{\em{)}} For an arbitrary Abelian surface $X$, 
the order $r$ of a finite automorphism on $X$ satisfies $\varphi (r)\leq 4$, 
which means that $r=1,2,3,4,5,6,8,10$. 
\end{prop}

Under Conjecture (\ref{conj}), we define the 
following constants. 

\begin{defn}\label{def3.3}
Assuming that $(F_l)$ holds true, we define 
$$
C_l:=\{c\in {\mathbb N}\, |\, c=|\rho_1(\xBir (S,0))|\}, 
$$
$$
D_l:=\{d\in {\mathbb N}\, |\, d=\xlcm(2,c),\ \text{where} \ c\in C_l\}, 
$$
$$
I_l:=\{e\in {\mathbb N}\, |\, e\  \text{is a divisor of } d\in D_l\}.  
$$
Assuming that $(F'_l)$ holds true, we define 
$$
I'_l:=\{\,c\in {\mathbb N}\, |\, c=|\rho_1 (g)| \}.
$$
\end{defn}

\begin{rem}\label{tuika}
The conjecture $(F_l)$ implies $(F'_l)$, and $(F_l)$ holds true for 
$l\leq 1$ by Proposition (\ref {prop20}) (1). 
In particular, we have that 
$$
I_0=\{1,2\}, \ \ I'_0=\{1\}, \ \ I_1=I'_1=\{1,2,3,4,6\}.   
$$
\end{rem}

\begin{prop}\label{surface}
The conjecture $(F_2)$ holds true with $B_2=66$ and 
$$
I'_2=\{r\in {\mathbb N}\,|\, \varphi (r)\leq 20, r\ne 60\}. 
$$
\end{prop}
\proof
Let $(S,0)$ be a normal surface with only canonical singularities such 
that $K_S\sim 0$. 
Let $f:T\to S$ be the minimal resolution. 
Note that $f$ is crepant, that is, $K_T=f^{*}K_S$. 
If $g\in \xBir (S,0)$, then $g':=f^{-1}\circ g\circ f\in 
\xBir (T,0)$. 
The discrepancy of every exceptional divisor over $T$ is positive and 
that of another divisor is non-positive. 
Since {\bb} map $g'$ does not change discrepancies, we have that 
$g'\in \xAut (T,0)=\xAut (T)$. 
By the classification of surfaces (see, for example, 
\cite[Theorem VIII.2]{Be}), 
$T$ is Abelian or $K3$. 
If $T$ is $K3$, then the second Betti number $b_2(T)=22$. 
If $T$ is Abelian, then the second Betti number $b_2(T)=6$. 
Therefore, $(F_2)$ holds true with $B_2=66$ 
by the proof of \cite [Proposition 14.4]{Ue} and 
$I'_2=\{r\in {\mathbb N}\,|\, \varphi (r)\leq 20, r\ne 60\}$ 
by Proposition (\ref {prop20}) (2). 
\endproof

\section{Indices of lc pairs with standard coefficients}\label{se4}

In this section, we use the following notations and the log MMP in dimension 
$n$ freely. All the results are stated for algebraic varieties. 
For analytic spaces, we recommend the reader to see Section \ref{se6}.   

\begin{notation}[{cf.~\cite[\S 2]{Sh}}]\label{notation2} 
Let $(P\in X,\Delta)$ be an $n$-dimensional 
log canonical pair such that $\Delta=\sum 
d_i\Delta_i$ has only standard coefficients. 
From now on, we may shrink $X$ around $P$ without mentioning it. 
If $\xInd (P\in X,\Delta)=a$, then $a(K_X+\Delta)\sim0$. 
The corresponding finite cyclic cover 
$$
\pi:(Q\in Y,D)\to (P\in X,\Delta)
$$ 
of degree $a$ is ramified only over the components of $\Delta_i$ of 
$\Delta$ with $d_i<1$ (see \cite[2.3, 2.4]{Sh}). 
Since $\Delta$ has only standard 
coefficients, $D$ is reduced and $D=\pi^{*}\llcorner \Delta \lrcorner$. 
We say that this $\pi:(Q\in Y,D)\to (P\in X,\Delta)$ is the 
{\em{index $1$ cover}} of the log divisor $K_X+\Delta$. 
By the construction, $K_Y+D=\pi^{*}(K_X+\Delta)$ has index $1$ and 
the index $1$ cover is unique up to \'etale equivalences. 
Let $G$ be the cyclic group associated to the cyclic cover 
$\pi:Y\to X$. 
Then we have the followings: 
$$
(Q\in Y,D)/G\simeq (P\in X,\Delta), 
$$
\begin{gather}
\tag{$\diamondsuit$}
(\pi_*{\mathcal O}_{Y}(m(K_Y+D)))^{G}\simeq {\mathcal O}_{X}(mK_X+
\llcorner m\Delta\lrcorner),
\end{gather}
where $m\in {\mathbb {Z}}_{\geq 0}$. 
From now on, we assume that $P\in \xCLC(X,\Delta)$. 
By \cite[Chapter 20]{FA}, \cite [\S 2]{Sh}, or 
\cite[Proposition 5.20]{KM}, $(Q\in Y,D)$ is log canonical 
and $P\in \xCLC(X,\Delta)$ is 
equivalent to $Q\in \xCLC(Y,D)$. Therefore, we may assume that 
$\xLLCC(Y,D)=Q$ without loss of generality. 
\end{notation}

\begin{prop}\label{fuan}
Let $m_0$ be a non-negative integer. 
Let $s$ be a $G$-invariant generator of ${\mathcal O}_{Y}(m_0(K_Y+D))$. 
Then $m_0(K_X+\Delta)$ is Cartier. 
In particular, $m_0\Delta$ is integral. 
\end{prop}
\proof
By ($\diamondsuit$), there exists a generator $t$ of 
${\mathcal O}_{X}(m_0 K_X+\llcorner m_0\Delta\lrcorner)$ such that 
$\pi^{*}t=s$. 
In particular, $m_0 K_X+\llcorner m_0\Delta\lrcorner$ is 
Cartier. 
Let $l$ be a sufficiently divisible positive integer such that 
$l\Delta$ is integral. 
Since $t$ is a generator of  
${\mathcal O}_{X}(m_0 K_X+\llcorner m_0\Delta\lrcorner)$, 
$t^{\otimes l}$ is that of 
${\mathcal O}_{X}(lm_0 K_X+l\llcorner m_0\Delta\lrcorner)$. 
On the other hand, $s^{\otimes l}$ is a generator of 
${\mathcal O}_{Y}(lm_0(K_Y+D))$ and 
$(\pi_*{\mathcal O}_{Y}(lm_0(K_Y+D)))^{G}\simeq {\mathcal O}_{X}(lm_0(K_X+
\Delta))$. 
Therefore, $t^{\otimes l}$ is also a generator of 
${\mathcal O}_{X}(lm_0(K_X+\Delta))$. Thus 
$\llcorner m_0 \Delta\lrcorner=m_0\Delta$, and $m_0(K_X+\Delta)$ is Cartier.  
\endproof

\begin{prop}[{cf.~\cite [Lemma 3.3]{I}}]\label{123}
Let $m_{Q}$ be the maximal ideal of $Q$ and 
$\rho:G\to \xGL(\omega_{Y}(D)/m_{Q}\omega_{Y}(D))$ 
the canonical representation. 
Let $m_0:=|\rho(G)|$. 
Then $\xInd (P\in X,\Delta)=m_0$. 
\end{prop}
\proof
The cyclic group $G$ acts on $\omega_{Y}(D)^{\otimes m_0}/m_{Q}\omega_{Y}(D)
^{\otimes m_0}$ trivially. 
Let $s'$ be a generator of $\omega_{Y}(D)^{\otimes m_0}/m_{Q}\omega_{Y}(D)
^{\otimes m_0}$ and $s''$ a lift of $s'$ in $\omega_{Y}(D)
^{\otimes m_0}$. 
Put 
$$
s:=\frac {1}{|G|}\sum _{\sigma \in G}\sigma^{*}s''.
$$ 
Then $s$ is a lift of $s'$ and a $G$-invariant generator of 
${\mathcal O}_{Y}(m_0(K_Y+D))$. 
By Proposition (\ref {fuan}), 
we obtain that $r:=\xInd (P\in X,\Delta)$ divides $m_0$. 
On the other hand, by considering the pull-back of the 
generator of $r(K_X+\Delta)$, we obtain that $m_0$ divides $r$. 
So we get $\xInd (P\in X,\Delta)=m_0$. 
\endproof

The following lemma is a special case of Lemma (\ref{taisetu2}). 

\begin{lem}\label{taisetu}
Let $h:(V,F)\to (Q\in Y,D)$ be a very good resolution. 
Let 
$$
V = V^0 \stackrel{p_1}{\dashrightarrow} V^1
 \stackrel{p_2}{\dashrightarrow} \cdots
 \stackrel{p_i}{\dashrightarrow} V^i
 \stackrel{p_{i+1}}{\dashrightarrow} V^{i+1}
 \stackrel{p_{i+2}}{\dashrightarrow}
 \cdots \stackrel{p_{l-1}}{\dashrightarrow}
 V^{l-1} \stackrel{p_l}{\dashrightarrow} V^l = Z
$$
be the $(K_V+F^{B})$-log MMP over $Y$.  
We denote $F_{0}^{B}=F^{B}, F_{0}=F$, and 
$F_{i}^{B}=p_{i*}F_{i-1}^{B}, 
F_{i}=p_{i*}F_{i-1}$, for every $i$ and $F_{l}=E$. 
Then $f:(Z,E)\to (Y,D)$ is a $\bQ$-factorial dlt model 
and $g:=p_{l}\circ p_{l-1}\circ\cdots\circ p_{1}$ induces an isomorphism 
at every generic point of center of log 
canonical singularities for the pair $(V,F)$. 
We note that $E=F_{l}=F_{l}^{B}$. 
Furthermore, $\xLLCC(Z,E)=E^{c}$. 
\end{lem} 
\proof 
Since $K_Y+D$ is Cartier, the effective part of $F$ is $F^{B}$. 
Therefore, by Lemma (\ref{taisetu2}), we get the result. 
\endproof

The next proposition is very important. 
We prove it by using Proposition (\ref{keyprop}). 
Note that, if $(Q\in Y,0)$ is an isolated singularity, then 
this proposition is obvious by Lemma (\ref{lem3}). 

\begin{prop}\label{prop4.2}
Let $h:(V,F)\to (Q\in Y,D)$ be a very good resolution. 
Then $\xLLCC (V,F)=F^{c}$ is connected.
\end{prop}

\proof
Since $h$ is a very good resolution, we have $\xLLCC (V,F)=F^{c}$. 
Run the $(K_V+F^{B})$-log MMP over $Y$. 
We get a $\bQ$-factorial dlt model 
$f:(Z,E)\to (Y,D)$ (see Lemma (\ref{taisetu})). 
We put $(T,\Xi):=(V,F)$, $(S,\Theta):=(Z,E)$, 
and $(o\in R):=(Q\in Y)$ and apply Proposition 
(\ref{keyprop}). 
The conditions (1), (2) and (5) in Proposition (\ref{keyprop}) are 
satisfied since $K_Y+D$ is Cartier and 
$K_V+F=h^{*}(K_Y+D)$, $K_Z+E=f^{*}(K_Y+D)$ and $Q\in \xCLC(Y,D)$. 
The condition (3) is in the definition of the very good resolution and (4) 
has been already checked in Lemma (\ref{taisetu}). 
Therefore, we can apply Proposition (\ref {keyprop}). 
Thus we obtain that $F^{c}$ is connected. 
\endproof

The following is a corollary of Proposition (\ref{prop4.2}). 
However, we don't use it for the proof of the main result. 

\begin{cor}\label{0}
Let $h':(V',F')\to (P\in X,\Delta)$ be a very good resolution. 
Then $\xLLCC(V',F')={F'}^{c}$ is connected. 
\end{cor}
\proof
Since $h'$ is a very good resolution, we have that $\xLLCC(V',F')={F'}^{c}$. 
Let $h:(V,F)\to (Q\in Y,D)$ be a very good resolution which 
factors $Y\times_{X}V'$. By Proposition (\ref{prop4.2}), $\xLLCC(V,F)
=F^{c}$ is 
connected. 
Since $h''(F^{c})={F'}^{c}$, where $h'':V\to V'$, 
we obtain that ${F'}^{c}$ is connected. 
\endproof

\begin{prop}\label{prop4.3}
There exists a $\bQ$-factorial 
dlt model $f:(Z,E)\to (Y,D)$, that is, $(Z,E)$ is dlt and 
$$
K_Z+E=f^{*}(K_Y+D),
$$ 
such that the following conditions are satisfied: 
\begin{enumerate}
\item[(1)] $\xLLCC (Z,E)=E^{c}$, 
\item[(2)] $E^{c}$ is connected and $(E^{c},(E-E^{c})|_{E^{c}})$ is sdlt, 
\item[(3)] $K_Z+(E-E^{c})|_{E^{c}}=(K_Z+E)|_{E^{c}}=f^{*}(K_Y+D)
|_{E^{c}}\sim 0$.  
\end{enumerate}
\end{prop}

\proof 
Let $f:(Z,E)\to (Y,D)$ be a $\bQ$-factorial dlt model constructed in 
Lemma (\ref{123}). 
Let $h':(U,H)\to (Z,E)\to (Q\in Y,D)$ be a very good resolution. 
Then $H^{c}$ is connected by Proposition (\ref {prop4.2}) and 
$g'(H^{c})=E^{c}$, where $g':U\to Z$. 
Therefore, $E^{c}$ is connected. 
Since $(Z,E)$ is $\bQ$-factorial and dlt, $E^{c}$ is Cohen-Macaulay and 
$(E^{c},(E-E^{c})|_{E^{c}})$ is sdlt by the adjunction. 
We note that $\xDiff (E-E^{c})=(E-E^{c})|_{E^{c}}$ since 
$(Z,E)$ is dlt and $K_Z+E$ is Cartier (see, for example, \cite[16.6 
Proposition]{FA}).  
\endproof

\begin{rem}\label{rrr}
By using Lemma (\ref{taisetu2}) and Corollary (\ref{0}), we get a 
similar result about $(P\in X,\Delta)$ by the same proof as that of 
Proposition (\ref{prop4.3}). 
That is, there exists a $\bQ$-factorial dlt model 
$f':(Z',E')\to (X,\Delta)$, that is, $(Z',E')$ is dlt and 
$K_{Z'}+E'={f'}^{*}(K_X+\Delta)$, 
such that the following conditions are satisfied: 
\begin{enumerate}
\item[(1)] $\xLLCC (Z',E')={E'}^{c}$, 
\item[(2)] ${E'}^{c}$ is connected and $({E'}^{c},
\xDiff (E'-{E'}^{c}))$ is sdlt, 
\item[(3)] $K_{Z'}+\xDiff (E'-{E'}^{c})=(K_{Z'}+E')|_{{E'}^{c}}
={f'}^{*}(K_X+\Delta)|_{{E'}^{c}}\sim_{\bQ} 0$.  
\end{enumerate}
See also Lemma (\ref{11}).
\end{rem}

\begin{lem}\label{lem14}
We have the following isomorphisms: 
\begin{eqnarray*}
{\omega_Y(D)}^{\otimes m}/m_{Q}{\omega_Y(D)}^
{\otimes m}&\simeq& H^{0}(E^{c}, f^{*}m(K_Y+D)|_{E^{c}})\\ 
&\simeq& 
H^{0}(E^{c}, m(K_{E^{c}}+(E-E^{c})|_{E^{c}}))\\&\simeq& 
H^{0}(E^{c}, {\mathcal O}_{E^{c}})\simeq {\mathbb C}, 
\end{eqnarray*}
where $m_Q$ is the maximal ideal of $Q$. 
\end{lem}

\proof 
We consider the following exact sequence, 
\begin{gather}
\tag{$\spadesuit$} 0\to {\mathcal O}_{Z} (-E^{c})\to 
{\mathcal O}_{Z}\to 
{\mathcal O}_{E^{c}}\to 0.
\end{gather}
Thus we obtain 
$$
{\mathcal O}_{Y}/ f_{*}{\mathcal O}_{Z} (-E^{c})\simeq {\mathbb C}.
$$
Note that $E^{c}$ is connected and $f(E^{c})=Q$. 
Since it is obvious that 
$m_Q\subset f_{*}{\mathcal O}_{Z} (-E^{c})$, so we obtain 
$m_Q=f_{*}{\mathcal O}_{Z} (-E^{c})$. 
By tensoring ${\mathcal O}_{Z}(f^{*}m(K_Y+D))$ to ($\spadesuit$) 
and taking direct images, 
we get 
\begin{align*}
0 \to  m_Q\otimes {\mathcal O}_{Y} (m(K_Y+D)) \to  
{\mathcal O}_{Y} (m(K_Y+D))\\  
 \to  
H^{0}(E^{c}, f^{*}m(K_Y+D)|_{E^{c}}) \to  \cdots .
\end{align*}
Therefore, we obtain the required isomorphisms. 
\endproof 

\begin{prop}\label{lem15}
Let $\sigma$ be an element of $G$. 
The {\bb} automorphism $\sigma:(Y,D)\to (Y,D)$ induces a {\bb} map 
$\sigma_{Z}:=f^{-1}\circ \sigma \circ f:(Z,E)\dashrightarrow (Z,E)$. 
Let $\alpha,\beta: (V,F)\to (Z,E)$ be a common resolution of 
$\sigma_{Z}:(Z,E)\dashrightarrow (Z,E)$ such that 
$f\circ \alpha, f\circ \beta: (V,F)\to (Q\in Y,D)$ are 
very good resolutions. 
Then we get the following commutative diagram: 

$$
\begin{CD}
H^{0}(E^{c},m(K_{E^{c}}+(E-E^{c})|_{E^{c}}))   @<
\text {$\sim$}< \text{$\sigma_{E^c}^*$}< 
H^{0}(E^{c},m(K_{E^{c}}+(E-E^{c})|_{E^{c}}))     \\
@A\text{$f^*$}A\text{$\simeq$}A      @A\text{$\simeq$}A\text{$f^*$} A  \\
\omega_{Y}(D)^{\otimes m}/m_Q\omega_{Y}(D)^{\otimes m}   @<\text {$\sim$}<
\text{$\sigma^*$}<  \omega_{Y}(D)^{\otimes m}/m_Q\omega_{Y}(D)^{\otimes m}.    
\end{CD}
$$
Here $\sigma _{E^c}^{*}=(\alpha^{*})^{-1}\circ \beta^{*}$, 
where $\alpha^*,\beta^*: H^{0}(E^{c},m(K_{E^{c}}+(E-E^{c})|_{E^{c}})) 
\simeq H^{0}(F^{c}, m(K_{F^{c}}+(F-F^{c})|_{F^{c}}))$. 
\end{prop}

\proof
We note the following isomorphisms: 
$$
H^{0}(F^{c},m(K_{F^{c}}+(F-F^{c})|_{F^{c}}))\simeq \mathbb C, 
$$
$$
H^{0}(E^{c},m(K_{E^{c}}+(E-E^{c})|_{E^{c}}))\simeq \mathbb C, 
$$ 
$$\omega_{Y}(D)^{\otimes m}/m_Q\omega_{Y}(D)^{\otimes m}\simeq \mathbb C. 
$$ 
Then we get the above commutative diagram by Lemma (\ref {lem14}). 
\endproof

\begin{prop}\label{prop4.10'}
In Proposition (\ref{lem15}), if there exists a non-zero {\ad} section 
in $H^{0}(E^{c},m_0(K_{E^{c}}+(E-E^{c})|_{E^{c}}))$, then $G$ acts on 
$\omega_{Y}(D)^{\otimes m_0}/m_Q\omega_{Y}(D)^{\otimes m_0}$ trivially. 
\end{prop}
\proof
This can be checked by the same argument as that of \cite [Lemma 4.7]{Fj}. 
Let $s\in H^{0}(E^{c},m(K_{E^{c}}+(E-E^{c})|_{E^{c}}))$ be a non-zero 
{\ad} section. 
It is sufficient to prove $\alpha^*s=\beta^*s$ in 
$H^{0}(F^{c},m_0(K_{F^{c}}+(F-F^{c})|_{F^{c}}))$. 
Let $F^{1}$ be any irreducible component of $F^{c}$. 
By applying \cite[Claim ($B_n$)]{Fj} repeatedly, we get $F^{2}\Subset F^{1}$ 
or $F^{2}=F^{1}$ such that $\alpha: F^{2}\to \alpha(F^{2})$ and 
$\beta:F^{2}\to \beta(F^{2})$ are {\bb} morphisms and 
$H^{0}(F^{1},m_0(K_V+F)|_{F^{1}})\simeq 
H^{0}(F^{2},m_0(K_V+F)|_{F^{2}})$ 
(see the proof of \cite[Lemma 4.7]{Fj}). 
Since $\alpha(F^{2})$ is $B$-birationally equivalent to $\beta(F^{2})$ 
and $s$ is a non-zero {\ad} section, 
we obtain that $\alpha^{*} (s|_{\alpha (F^{2})})=
\beta^{*}(s|_{\beta(F^{2})})$ in $H^{0}(F^{2}, m_0(K_V+ F)|_{F^{2}})$. 
Therefore, we have that $\alpha^{*}s=\beta^{*}s$ in 
$H^{0}(F^{1}, m_0(K_V+ F)|_{F^{1}})$. 
Thus we obtain that $\alpha^{*}s=\beta^{*}s$ in 
$H^{0}(F^{c}, m_0(K_{F^{c}}+ (F-F^{c})|_{F^{c}}))$. 
\endproof

\begin{defn}\label{defofmu}
Let $(P\in X,\Delta)$ and $(Q\in Y,D)$ be as in Notation (\ref{notation2}). 
Let $h:(U,H)\to (Q\in Y,D)$ be a log resolution such that 
$K_U+H=h^{*}(K_Y+D)$. 
We define 
$$
\mu=\mu (P\in X,\Delta):= \min \{\dim W\, |\, W\in \xCLCC(U,H)\}.
$$ 
It is obvious that $0\leq\mu \leq \dim X-1$. 
We note that the index $1$ cover $(Q\in Y,D)$ is defined uniquely up to 
\'etale equivalences. 
By Lemma (\ref{99}) (see also Section \ref{se6}), 
$\mu$ is independent of the choice of the resolution. 
Therefore, $\mu (P\in X,\Delta)$ is well-defined. 
\end{defn}

\begin{rem}\label{remdual}
When $(P\in X,0)$ is an $n$-dimensional isolated log canonical 
singularity, which is not log terminal, 
Shihoko Ishii defined the lc singularity of type $(0,i)$ by using the 
mixed Hodge structure of the simple normal crossing 
variety $H^{B}$ (see \cite [2.7]{I}), 
where $h:(U,H)\to (Q\in Y,0)$ is a log resolution 
as in Definition (\ref {defofmu}). 
She also proved that 
$$
\dim \Gamma _{H^{B}} =n-1-i,
$$ 
where $\Gamma _{H^{B}}$ is the dual graph of $H^{B}$ 
by \cite [Theorem 2]{I1} and \cite{I4}. 
By the definition of $\mu (P\in X,0)$, 
we have 
$$
\dim \Gamma _{H^{B}} =n-1-\mu.
$$ 
Therefore, we get $\mu =i$ when $(P\in X,0)$ is an isolated 
log canonical singularity. 
Furthermore, if the log MMP holds true in dimension $\leq n$, 
then the above dual graph $\Gamma _{H^{B}}$ is pure 
$(n-1-\mu)$-dimension by Lemma (\ref{99}) 
and Lemma (\ref{minimal}) below. 
\end{rem}

\begin{prop}\label{prop4.8}
Assume that $\mu (P\in X,\Delta)\leq n-2$. 
If $(F_{\mu})$ holds true, then there is a non-zero {\ad} section 
$s\in H^{0}(E^{c},m_0(K_{E^{c}}+(E-E^{c})|_{E^{c}}))$ with 
$m_0\in D_{\mu}$. 
In particular, $s$ is $G$-invariant. 
Thus, $\xInd (P\in X,\Delta)\in I_{\mu}$.  
\end{prop}
\proof
Let $W$ be a compact minimal center of lc singularities for the pair 
$(Z,E)$. 
By the definition of  $\mu (P\in X,\Delta)$, $\dim W=\mu$. 
We note that $K_{E^{c}}+(E-E^{c})|_{E^{c}}\sim 0$. 
By using the adjunction repeatedly, we have that 
$(K_{E^{c}}+(E-E^{c})|_{E^{c}})|_{W}=K_W\sim 0$. 
So $(W,0)$ has only canonical singularities. 
By Lemma (\ref{minimal}) below, 
all the minimal compact centers of lc singularities are 
$B$-birationally equivalent to $(W,0)$. Therefore, 
$H^{0}(\amalg W,m_1(K_Z+E)|_{\amalg W})$ has a non-zero 
{\ad} section with $m_1:=|\rho_1(\xBir(W,0))|\in C_{\mu}$, 
where the sum runs over all the compact minimal centers of lc 
singularities for the pair $(Z,E)$. 
By applying Proposition (\ref{prop4.8.1}) below repeatedly, 
we get a non-zero {\ad} section with $m_0\in D_{\mu}$, 
where $m_0=\xlcm (2,m_1)$. 
Therefore, by Propositions (\ref{fuan}), (\ref{lem15}), (\ref{prop4.10'}) and 
the proof of Proposition (\ref{123}), we get $\xInd (P\in X,\Delta)$ 
is a divisor of $m_0$. 
In particular, $\xInd (P\in X,\Delta)\in I_{\mu}$.    
\endproof

We note that, in Lemmas (\ref{minimal}) and (\ref{11}) and Proposition (\ref
{prop4.8.1}), 
we use Proposition (\ref {conn}) and \cite[Proposition 2.1]{Fj}, 
which need the log MMP. 

\begin{lem}\label{minimal}
Let $(S,\Delta)$ be a proper connected sdlt $n$-fold with 
$K_S+\Delta\sim 0$. 
Let $\mu :(S',\Theta)=\amalg(S_i,\Theta_i)\to (S,\Delta)$ 
be the normalization. 
Then all the minimal centers of log canonical singularities for the pair 
$(S',\Theta)$ have the same dimension and are $B$-birationally 
equivalent to each other. 
\end{lem}

\proof
We prove this lemma by induction on $n$. 
If $n=1$, then this is trivial. 
Let 
$$
d:=\min\{\dim W\,|\,W\in \xCLC (S',\Theta)\}. 
$$

If $d\leq n-2$, then $(S_i,\Theta_i)$ is in 
Case (1) in Proposition (\ref{conn}) 
and $\Theta_i$ is not irreducible for every $i$. 
By applying the induction to $(n-1)$-dimensional 
sdlt pair $(\Theta _i,0)$, we get the result. 

If $d=n-1$, then $(S_i,\Theta_i)$ is plt for every $i$. 
Therefore, all the minimal elements in $\xCLC(S', \Theta)$ have dimension 
$n-1$ and $B$-birationally equivalent to each other by Proposition 
(\ref {conn}). 
\endproof

For the main result, the above lemma is sufficient. 
However, in the above lemma, the assumption that $K_S+\Delta\sim 0$ can be 
replaced by  $K_S+\Delta\equiv 0$. 

\begin{lem}\label{11}
Let $(S,\Delta)$ be a proper connected sdlt $n$-fold with 
$K_S+\Delta\equiv 0$. 
Let $\mu :(S',\Theta)=\amalg(S_i,\Theta_i)\to (S,\Delta)$ 
be the normalization. 
Then all the minimal centers of log canonical singularities for the pair 
$(S',\Theta)$ have the same dimension and are $B$-birationally 
equivalent to each other. 
\end{lem} 
\proof
We prove this lemma by induction on $n$. 
If $n=1$, then this is trivial. 
Apply the inductive hyposesis to each connected component 
of $(n-1)$-dimensional sdlt pair $(\llcorner \Theta_i\lrcorner, 
\xDiff (\Theta_i-\llcorner \Theta _i\lrcorner))$ 
for every $i$ and 
use Lemma (\ref{99}) and \cite [Proposition 2.1 (0) (0.2)]{Fj}. 
Then we get the result. 
\endproof

In the proof of Proposition (\ref{prop4.8}), we used the following 
result, which is a special case of \cite[Proposition 4.3]{Fj}. 
In our situation, all the dlt pairs have Kodaira dimension $0$. 
So, \cite[Proposition 4.3]{Fj} can be modified as follows: 

\begin{prop}\label{prop4.8.1}
Let $(S,\Theta)$ be a proper dlt $n$-fold such that $K_S+\Theta\sim 0$. 
Assume that there exists a non-zero {\ad} section 
$u \in H^{0}(\Theta, m(K_S+\Theta)|_{\Theta})$. 
If $m$ is even, then we can extend $u$ to $v\in H^{0}(S, m(K_S+\Theta))$, 
that is, $v|_{\Theta}=u$. 
In particular, $v$ is a non-zero {\ad} section. 
\end{prop}
\proof 
This is a special case of Cases (1) and (4) in the 
proof of \cite[Proposition 4.3]{Fj} (for Case (4), see also 
\cite [Proposition 4.5]{Fj3}). 
For the latter part, see the case (3) in the proof of 
\cite [Theorem 3.5]{Fj}. 
\endproof 

In the case where $\mu (P\in X,\Delta)=
n-1$, we can prove a slightly stronger result 
by using the canonical desingularization theorem. 

\begin{prop}\label{prop4.9}
Assume that $\mu (P\in X,\Delta)=n-1$. 
Then there exists a projective birational morphism 
$f:(Z,E)\to (Y,D)$ from a dlt pair $(Z,E)$, which satisfies 
the following conditions.  
\begin{enumerate}
\item[(1)] $K_Z+E=f^{*}(K_Y+D)$.
\item[(2)] Let $\sigma \in G$. 
The {\bb} automorphism $\sigma:(Y,D)\to (Y,D)$ induces the {\bb} automorphism 
$\sigma_Z:=f^{-1} \circ \sigma \circ f 
:(Z,E)\to(Z,E)$ over $Y$ and the automorphism 
$\sigma _E:E\to E$.  
\item[(3)] The morphism $f$ is $G$-equivariant.
\item[(4)] There exists the following commutative diagram: 
$$
\begin{CD}
H^{0}(E,mK_{E})   @<
\text {$\sim$}< \text{$\sigma_{E}^*$}< 
H^{0}(E,mK_{E})     \\
@A\text{$f^*$}A\text{$\simeq$}A      @A\text{$\simeq$}A\text{$f^*$} A  \\
\omega_{Y}(D)^{\otimes m}/m_Q\omega_{Y}(D)^{\otimes m}   @<\text {$\sim$}<
\text{$\sigma^*$}<  \omega_{Y}(D)^{\otimes m}/m_Q\omega_{Y}(D)^{\otimes m}.    
\end{CD}
$$
\end{enumerate} 
\end{prop}
\proof
We take the canonical desingularization 
$h:(V,F)\to (Q\in Y,D)$ (see \cite{BM}). 
Since $h$ is canonical, $G$ acts on $V$. 
Since $\sigma_{V}:=h^{-1} \circ \sigma \circ h$ 
is a {\bb} automorphism, $G$ also acts on $F^{B}$. 
Run the $G$-equivariant log MMP with respect to $K_V+F^{B}$ over $Y$ 
(see \cite[Example 2.21]{KM}). 
Then we get a $G$-equivariant dlt model 
$f:(Z,E)\to (Y,D)$, 
that is, $G$ acts on $Z$ and $E$. 
Since $\mu (P\in X,\Delta)=n-1$, $E=E^{c}$ is irreducible. 
By Lemma (\ref {lem14}), we get (4). 
\endproof

\begin{prop}\label{prop4.10}
Let the notations be as in Proposition (\ref{prop4.9}). 
If the conjecture $(F'_{n-1})$ holds true, then 
there exists a non-zero $G$-invariant section $s\in H^{0}(E,m_0K_{E})$ 
with $m_0:=|\rho_1(G)|\in I'_{n-1}$. 
In particular, 
$G$ acts trivially on 
$\omega_{Y}(D)^{\otimes m_0}/m_Q\omega_{Y}(D)^{\otimes m_0}$. 
Thus we get $\xInd (P\in X,\Delta)=m_0\in I'_{n-1}$.  
\end{prop}
\proof
It is obvious by Propositions (\ref{prop4.9}), (\ref{fuan}), and 
(\ref{123}). 
\endproof

The following is the main theorem of this paper, which is a 
consequence of Propositions (\ref{lem15}), (\ref{prop4.10'}), (\ref{prop4.8}), 
(\ref{prop4.9}), and (\ref{prop4.10}). 

\begin{thm}\label{mainthm}
Assume the log MMP for dimension $\leq n$. 
Let $(P\in X,\Delta)$ be 
an $n$-dimensional lc pair such that 
$\Delta$ has only standard coefficients and $P\in \xCLC(X,\Delta)$. 
When $\mu(P\in X,\Delta)\leq n-2$ 
{\em{(}}resp.~$\mu (P\in X,\Delta)=n-1${\em{)}}, 
we assume 
$(F_{\mu})$ {\em{(}}resp.~$(F'_{n-1})${\em{)}} holds true. 
Then 
$$
\begin{cases}
\text{$\xInd (P\in X,\Delta)\in I'_{n-1}$}\\
\text{$\xInd (P\in X,\Delta)\in I_{\mu}$}\\
\end{cases}
\quad {\text{if}}
 \quad
\begin{cases}
\mu (P\in X,\Delta)=n-1, \\
\mu (P\in X,\Delta)\leq n-2.\\
\end{cases}
$$
\end{thm}

For three dimensional log canonical pairs, we 
obtain the following result as a corollary of Theorem (\ref{mainthm}) 
(for related results, see \cite {I} and \cite[1.10 Corollary]{Sh2}).  

\begin{cor}\label{maincor}
Let $(P\in X,\Delta)$ be a three dimensional lc pair such 
that $\Delta$ has only standard coefficients and $P\in \xCLC(X,\Delta)$.
Then 
$$
\begin{cases}
\text{$\xInd (P\in X,\Delta)\in \{1,2\}$}\\
\text{$\xInd (P\in X,\Delta)\in \{1,2,3,4,6\}$}\\
\text{$\xInd (P\in X,\Delta)\in I'_{2}$}\\
\end{cases}
\quad {\text{if}}
 \quad
\begin{cases}
\mu (P\in X,\Delta)=0, \\
\mu (P\in X,\Delta)=1,\\
\mu (P\in X,\Delta)=2,\\
\end{cases}
$$
where $I'_{2}=\{r\in {\mathbb N}\,|\, \varphi(r)\leq 20, r\ne 60\}$. 
In particular, if there exists 
$W\in \xCLC(X,\Delta)$ such that 
$P\subsetneq W$, then $\xInd(P\in X,\Delta)\in \{1,2,3,4,6\}$. 
\end{cor}

\begin{rem}\label{rem22}
Shihoko Ishii proved that for every $r\in {I'_2}=
\{r\in {\mathbb N}\,|\, \varphi(r)\leq 20, r\ne 60\}$, 
there exist three dimensional isolated log canonical 
singularities such that $\mu (P\in X,0)=2$ and 
$\xInd (P\in X,0)=r$ (see \cite [Theorem 4.15] {I}). 
For the singularities which satisfy $\mu \leq 1$, see Example (\ref{extuika}). 
\end{rem}

For two dimensional log canonical pairs, 
the following corollary seems to be well-known to specialists. 

\begin{cor}\label{nizigen}
Let $(P\in X,\Delta)$ be a two dimensional lc pair such 
that $\Delta$ has only standard coefficients and $P\in \xCLC(X,\Delta)$.
Then 
$$
\begin{cases}
\text{$\xInd (P\in X,\Delta)\in \{1,2\}$}\\
\text{$\xInd (P\in X,\Delta)\in \{1,2,3,4,6\}$}\\
\end{cases}
\quad {\text{if}}
 \quad
\begin{cases}
\mu (P\in X,\Delta)=0, \\
\mu (P\in X,\Delta)=1.\\
\end{cases}
$$
\end{cor}

\section{Examples}\label{se5}

In this section, we treat some examples of log canonical pairs. 

\begin{exmp}\label{ex7}
Let $X=(x^3+y^3+z^3=0)\subset {\mathbb C}^{4}=
\Spec\,{\mathbb C}[x,y,z,t]$. 
The cyclic group ${\mathbb Z}_{m}$ acts on $X$ as follows: 
$$
(x,y,z,t)\mapsto (\varepsilon x, \varepsilon y, 
\varepsilon z, \varepsilon t), 
$$ 
where $\varepsilon$ is a primitive $m$-th root of unity. 
Let $(o\in Y_m)$ be the quotient $X/{{\mathbb Z}_{m}}$, 
where $o$ is the image of $(0,0,0,0)\in X$. 
The cyclic group ${\mathbb Z}_{m}$ acts on 
$\omega_{X}={\mathcal O}_{X}(K_X)$ 
as follows: 
$$
\omega \mapsto \varepsilon \omega,
$$ 
where 
$$
\omega=\xRes \frac {dx\wedge dy\wedge dz\wedge dt}{x^3+y^3+z^3},
$$ 
which is a generator of $\omega_{X}$. 
Therefore we obtain that $\xInd (o\in Y_{m})=m$. 
This shows that the indices are not bounded if 
$o\notin \xCLC(Y_m,0)$. 
\end{exmp}

\begin{exmp}\label{Ex}
Let $X:=(x^2+y^3+z^7+t^6z^6=0)\subset {\mathbb C}^4$ 
and $o=(0,0,0,0)\subsetneq W:=\{(x,y,z,t)\in X\, |\, x=y=z=0\}$. 
Let $g:Z\to {\mathbb {C}}^{3}=\Spec\,{\mathbb C}[x,y,z]$ be the 
weighted blowing up at $(0,0,0)$ with the weight 
$(\wt x, \wt y, \wt z)=(3,2,1)$. 
Let $h:=g\times 1: Z\times {\mathbb{C}}\to  
{\mathbb {C}}^{4}=\Spec\,{\mathbb C}[x,y,z,t]$ 
and $Y$ the strict transform of $X$ by $h$. 
Let $f:=h|_{Y}:Y\to X$ and 
$E$ be the exceptional divisor of $f$. 
Then $K_Y=f^{*}K_X-E$ and $Y$ is smooth, and $(E,0)$ 
has only one lc singularity in $f^{-1}(o)$. 
So, by \cite[17.2 Theorem]{FA}, we obtain 
that $X$ is lc and $o$ and $W$ are centers of log canonical 
singularities for the pair $(X,0)$. 
\end{exmp}

\begin{exmp}\label{ex600}
Let $X={\mathbb {C}}^{4}$ and 
$$
\Delta:=\frac{1}{2}\Delta_1+\frac{2}{3}\Delta_2+\frac{7}{8}\Delta_3+
\frac{24}{25}\Delta_4+\frac{599}{600}\Delta_5, 
$$
where $\Delta_i$ is a general hyperplane and 
$o:=(0,0,0,0)\in \Delta_i$ for every $i$. 
Then $(o\in X,\Delta)$ is a four dimensional lc pair such that 
$\xLLC(X,\Delta)=o$ and $\xInd (o\in X,\Delta)=600$. 
\end{exmp}

\begin{exmp}\label{extuika}
Let $(P\in Z,0)$ be a two dimensional log canonical singularity which is 
not log terminal. 
Then, by Corollary (\ref{nizigen}), $\mu (P\in Z,0)=0$ or $1$, 
and $\xInd (P\in Z,0)=1,2,3,4$, or $6$. 
Let $X:=Z\times {\mathbb{C}}$ and $p:X\to {\mathbb{C}}$ be the 
second projection. 
Let $o\in {\mathbb{C}}$ and $H:=p^{*}[o]$, and $Q:=(P,o)\in X$. 
Then $(Q\in X,H)$ is log canonical and $Q\in \xCLC(X,H)$ by 
\cite [17.2 Theorem]{FA}. 
Furthermore, $\mu (P\in Z,0)=\mu (Q\in X,H)$ and $\xInd (P\in Z,0)=
\xInd (Q\in X,H)$. 
\end{exmp}

\section{Appendix}\label{se6}

In this appendix, we treat analytic germs of lc pairs. 
The main theorem holds true for analytic spaces. 

\begin{thm}[Analytic version of the main results]\label{analytic}
Theorem (\ref{mainthm}) and Corollaries (\ref{maincor}), 
(\ref{nizigen}) hold true 
even if $X$ is an analytic space. 
\end{thm} 

In Section \ref{se1}, we defined {\bb} map (resp.~morphism), 
$\xBir(X,D)$ and so on. 
For analytic spaces, {\bm} map (resp.~morphism), 
$\xBim(X,D)$ etc.~can be defined without difficulty. 
Details are left to the reader. 

For the relative log MMP over analytic germs, we recommend the reader to see 
\cite [Section 1]{Ka1} or \cite {N1}. 
By using this, Section \ref{se2} can be generalized to analytic spaces, 
which are projective over analytic germs. 
Note that in Propositions (\ref{prop4.3}) and (\ref{prop4.9}), $E^{c}$ is 
projective. 
So we don't have to modify Section \ref{se3}. 

In Section \ref{se4} and Lemma (\ref{99}), 
we often used \cite [Lemma 4.7]{Fj}. 
In \cite[Lemma 4.7]{Fj}, we used Szab\'o's resolution lemma (see \cite{Sz}). 
In \cite{Sz}, Szab\'o proved the resolution lemma for only algebraic 
varieties. 
Therefore, in the analytic case, we apply \cite[Theorem 12.4]{BM}, which 
contains the analytic version of Szab\'o's resolution lemma. 

In the proof of \cite [Lemma 4.7, Claim ($B_{n}$)]{Fj}, 
we cited \cite[Lemma 2.45]{KM}, which is an algebraic result. 
However, we can check that 
we didn't need it. 
In \cite[Claim ($B_n$)]{Fj}, we only treated {\bs} pairs. 
So, by using the following lemma, we can prove \cite[Claim ($B_n$)]{Fj} 
without using \cite[Lemma 2.45]{KM} 
(for related topics, see \cite [1.1]{Ku}, \cite {J}, and \cite 
[Chapter VI 1.4.7]{K}). 

\begin{lem}\label{saigo}
Let $(S,\Xi)\to (T,\Theta)$ be a {\bm} morphism between 
{\bs} pairs. 
Let $W$ be an irreducible component of $\Xi^{B}$. 
Then there exists a finite sequence of blowing-ups: 
$$
T^l \stackrel{p_l}{\longrightarrow} T^{l-1}
 \stackrel{p_{l-1}}{\longrightarrow} \cdots
 \stackrel{p_{k+1}}{\longrightarrow} T^k
 \stackrel{p_{k}}{\longrightarrow} T^{k-1}
 \stackrel{p_{k-1}}{\longrightarrow}
 \cdots \stackrel{p_{2}}{\longrightarrow}
 T^{1} \stackrel{p_1}{\longrightarrow} T^0 = T,
$$
whose centers are $W^{k}\in \xCLC(T^{k},\Theta^{k})$ and 
$W^l$ is a divisor. 
Here $p^{*}_k(K_{T^{k-1}}+\Theta^{k-1})=K_{T^{k}}+\Theta^{k}$ and $W^k\in 
\xCLC (T^k,\Theta^k)$ is a center associated to the divisor $W$ for every $k$. 
\end{lem}
\proof 
This can be checked easily by computing discrepancies by the same way 
as in \cite[Chapter VI 1.4.7]{K}. 
\endproof

Therefore, Theorem (\ref{analytic}) can be proved without difficulty. 
Details are left to the reader. 

\ifx\undefined\bysame
\newcommand{\bysame}{\leavevmode\hbox to3em{\hrulefill}\,}
\fi

\end{document}